\documentclass{article}\usepackage{amsfonts}
\usepackage{graphicx}
\begin{document}
\def \m {\vspace*{6pt}

\noindent}
\large\centerline{ Finite generation and the Gauss process \footnote
{John Atwell Moody, Maths, Warwick University;  
moody@maths.warwick.ac.uk }} 

\m \m {\narrower Abstract: Convergence of the Gauss resolution process 
for a complex singular foliation of dimension $r$ 
is shown to be equivalent to finite type of a graded sheaf  
which is built using base $(r+2)$ expansions
of integers. As  applications it is calculated which foliations 
coming from split semisimple representations of commutative Lie algebras
can be resolved  torically with respect to an eigenspace decomposition
and it is shown that Gaussian
resolutions stabilize for  irreducible projective varieties with
foliations of dimension $r$  for which $(r+1)H+K$
is a finitely-generated divisor of Iitaka dimension less than two where $H$ is a hyperplane
section  and $K$ a canonical divisor of the foliation. Another
application is that for normal irreducible   complex projective
varieties  with very ample divisor $H$ and a resolvable foliation, 
there are functorial locally closed conditions
on vector subspaces $X\subset |iH|$ which hold for large $i$
and ensure that 
 blowing up the base locus of $X$ and one further
Gaussian blowup  resolves the foliation. 

  \m}\m
\centerline{1. Introduction.}  
\m  Let $V$ be an irreducible algebraic variety over a field $k$
of characteristic zero, furnished with a singular foliation.
\m 
 The {\it Gauss blowup}
of the foliated variety $V$ is the image of the
 Gauss rational map to a Grassmannian, or rather the lowest domain
of definition of that map.
In reasonable cases a foliation which can be resolved can be resolved by
its Gauss map, but this is not always so.  A necessary and sufficient
condition for resolvability of a singular foliation 
by an arbitrary locally projective birational 
morphism 
in characteristic zero was described  in an earlier paper [1], that there is
an ideal sheaf $I$ and a number $N$ so that
$$I^N{\cal J}(I)^{r+2}= I^N{\cal J}(I{\cal J}(I)).$$
\m The motivating goal of this paper is to show that if the condition 
holds for some value of $N$ then  there is some possibly 
different choice of $I$ so that it holds with $N=0.$ 
\m Working, initially from affine examples, we shall construct a
sheaf $R(I)$ of graded algebras, depending on $I$, making use of base $(r+2)$ expansions
of integers, such
 that finite generation (`finite type')
of $R(I)$ over $k$  is equivalent
to convergence of the  resolution process where first one blows up $V$
along  $I$
and then follows by iterating the Gauss process. When this is the case,
the locally projective  morphism associated 
to the graded sheaf turns out
not to be the limit (=last stage) of the tower, but rather it is
a variety which converts the next-to-last stage to the last
stage by pulling back.  
\m It  follows from the construction  that
resolving the foliation by a locally projective morphism is  
equivalent to solving the simpler equation
${\cal J}(I)^{r+2}={\cal J}(I{\cal J}(I)).$ 
The simpler condition is now neither necessary
nor sufficient for $I$ itself to be a resolving ideal.
\m
 The operator ${\cal J}$ in [1],  which operated
upon fractional ideals,  depended upon a  basis
of derivations of the rational function field of $V$ over $k.$
We will replace it with 
a more functorial version of the same operator:
therefore with a functor $F.$ 
\m In this introduction, which will be light reading, we won't define the
functor $F;$ it is enough to know that it operates on torsion
free rank one coherent sheaves
 on $V. $ Then the graded sheaf $R$ 
also  depends functorially
on a torsion free coherent sheaf $I, $  and let us describe
how things fit together.
\m 
Each graded part  $R_i$ 
is  also to be  functor operating on torsion free rank one coherent
sheaves, such that the $i$'th term of $R(I)$ is 
the functor $R_i$ applied to the sheaf $I.$   Here is how the
functor $R_i$ is defined, starting 
 with the functor $F.$ If $H$ and $G$ are functors
acting on torsion free coherent sheaves on $V$ of rank one, we
will always define
the product $HG$ to be the most obvious functor which assigns to each $J$
the sheaf $H(J)G(J)$ by which we mean the tensor product
of $H(J)$ and $G(J)$ reduced modulo torsion. Note that the identity
functor $id$ is not  an identity element for the
multiplication. 
\m In terms of this product
operation,
define a sequence of functors $L_i$ by  the inductive rules
$$L_0=F$$
$$L_{i+1}=F(\ id\  L_0...L_i)$$
\m Then for each number $i$  define the functor
$$R_i=L_0^{a_0}L_1^{a_1}...L_s^{a_s}$$
where $i=a_0+a_1(r+2)+...+a_s(r+2)^s$ with $0\le a_\alpha<(r+2)$
is the base $r+2$ expansion of the degree $i.$ For each
choice of sheaf $I$ we show the $R_i(I)$ fit together to be a graded sheaf
of algebras.
\m  
The multiplication law
 is to be  based on the carrying operation which takes place 
when two finite geometric
series to the base of $(r+2)$ are added. The carrying 
operation in the $s$'th place, 
for a torsion free coherent sheaf $J,$ 
wll require  a map
$$F(JL_0(J)...L_s(J))^{r+2} \to F(JL_0(J)...L_{s+1}(J))$$
satisfying suitable compatibility conditions for the various values of $s$
coming from the ring axioms. The fact that the left and
right sides of the map above have the same degree in the
graded sheaf $R$  requires  
that the  action of $F$
on degrees will be be the affine transformation 
of multiplying by $r+1$ and then adding $1$,  corresponding to the
numerical identity
$$(1+(r+1)(1+(r+2)+...+(r+2)^s))(r+2)  \ \ \ \ \ \ \ \ \ 
$$ $$\ \ \ \ \ \ \ \ \ \ \ =1+(r+1)(1+(r+2)+...+(r+2)^{s+1}).$$
A goal of later sections will be to begin to understand
a manner in which  the graded sheaf $R(I)$ and the
numerical identity, may have a meaningful representation together in the context
of algebraic geometry.

\m  The more precise statement of our result,
now, is that starting with a singularly foliated irreducible
variety $V$ over a field $k$ of characteristic zero, 
furnished with  a torsion free rank one coherent sheaf $I,$ 
 the graded
sheaf $R(I)$ is finite type over $k$ 
if and only if that sequence of blowups of $V$
stabilizes, which consists of  blowing up $I$  and then
taking successive Gaussian blowups with respect to the foliation. 
\m As a test application of the picture, we will construct
the sheaf $R({\cal O}_V)$ when $V$ is the cusp, and also
in the case of the
simplest non-resolvable foliation of the plane and other
examples. Then we  will 
apply the rule to determine precisely which of the foliations
coming from split  semisimple representations of  commutative
Lie algebras can be resolved torically, with respect to  the
toric structure coming from an eigenspace decomposition.
The answer will show us
that those which can be  resolved torically  can all be resolved by one sing
le Gauss blowup.
\m Next, taking $V$ to be a normal irreducible quasiprojective variety, and
$I$ to be a very ample line bundle, we will show that $R(I)$ is
always generated by global sections.
Therefore the finite type property for $R(I),$ and  in turn also
convergence of the Gauss process itself, is equivalent
to  finite generation of the algebra of global sections.  
\m Next assume the quasiprojective variety $V$  is actually
projective and 
set $D=(r+1)H+K_V$ where $H$ is a hyperplane section and $K_V$
 is the canonical divisor of the foliated variety 
(in the sense of the introductory pages Bogomolov \& McQuillan [16], for example).
\m It is possible to practically dispense with
the cases when $D$ is a finitely generated divisor
of  Iitaka dimension
zero or 1. 
\m
If the Iitaka dimension of $D$ 
is equal to zero, then the graded algebra in question
is an integral domain of transcendency one degree over $k;$  
this can be nothing but a polynomial algebra $k[T]$ in one variable
$T$ if $k$ is algebraically closed. This is finite type over $k$
so the Gauss process for $V$ must stabilize; but analyzing
further we see that in this case the foliation is always nonsingular.
\m This is because the 
letter $T$ can be seen as a global section of ${\cal O}_V(K_V+(r+1)H),$
and the fact that the degree one monomials of the polynomial
algebra are spanned by $T$ means this global section spans a complete
linear system which generates the 
sheaf   of
 $r$'th  exterior power of the differentials of $V$
along the foliation
mod torsion.  The sheaf is therefore
trivial of rank one, and the foliation nonsingular.
\m If we assume in addition that the foliation is codimension
zero, ie that $r=dim(V),$ this would imply that the
nash blowups of $V$ itself stabilize, but
this is for a trivial reason, for then  $V$ can be nothing
but projective space (see [18] chapter 3). In our situation 
this is easily seen as follows, as
  Kobayashi and
Ochiai [6] applies  to the nonsingular variety $V.$
Namely, writing $c$ for the Chern class of $H$ the Chern
character of ${\cal O}_V(iH)$ is 
$$ch(iH)=1 + ic + i^2/2! c^2 + ... + i^r/r! c^r.$$
By Hirzebruch Riemann-Roch 
$$\chi(iH)=ch(iH)Td(V)=Td_n(V)+ i\cdot Td_{n-1}(V).c + ... +i^r/r!\cdot c^r.$$
Since $K=-(r+1)H,$ Kodaira vanishing shows that the
right side is equal to zero 
for $i= -1, ..., -r$  and 1 for $i=0$.
Then the rights side for all $i$ is
$$\chi({\cal O}_V(iH))={1\over {r!}}(i+r)(i+r-1)...(i+1)
=\left(\matrix{r+i\cr r}\right)$$
Then 
$$H^r=c^r=1$$
showing  
$H^r=1.$ This shows that $V$ meets an intersection of $r$
hyperplanes at a single point, so is a linear projective space.

\m
For cases when $D$ has nonzero Iitaka dimension, one
needs to relate 
the sheaf of graded algebras  ${ R}({\cal O}_V)$ is
to   the 
graded `pluricanonical sheaf' 
 $\oplus_i {\cal O}_V(iK_V)).$ 
\m 
Let then $V$ be a  normal  irreducible projective variety over a field $k$
of characteristic zero, with a singluar foliation, 
and $H$ a very ample divisor. 
Let ${ R}$ be the graded sheaf corresponding to the structure
sheaf ${\cal O}_V.$
Each ${ R}_i$ can be represented as a subsheaf
of ${\cal O}_V(iK_V),$  compatibly with multiplication
of sections, so ${ R}$ can be represented as a subsheaf of
graded rings of $\oplus_i {\cal O}_V(iK_V).$  Correspondingly the 
sheaf of rings
${ R} ({\cal O}_V(H))$ will be 
a subsheaf of $\oplus_i {\cal O}_V(iD)$ for 
$$D=K_V+(r+1)H$$ and we will find that the subsheaf  
is generated, as a graded
coherent sheaf, by global sections.
\m
In this way we obtain  a subalgebra of 
$$\oplus_{i=0}^\infty H^0(V, {\cal O}_V(iD))$$ for this divisor $D$
whose finite generation controls stability of the Gaussian resolution
process for the foliated variety $V.$ 
\m  We may take the subalgebra  to be
generated by particular
global sections whose degrees are powers of $(r+2),$ defined
explicitly, inside the full algebra
of global sections of $R({\cal O}_V(H)).$

\m Here is a geometric interpretation. 
Abbreviate $L_i=L_i({\cal O}_V(H)), $ and $R_i=R_i({\cal O}_V(H))$
omitting the sheaf ${\cal O}_V(H)$ by abuse of notation.
Since  $L_i=R_{(r+2)^i}$  each $L_i$ is generated
by global sections.
\m The product 
$$\Gamma(L_0)... \Gamma(L_s) \ \ \ \ \ \ \ \ \ $$
$$ \ \ \ \ \ \ \ \ \ \ \ \subset  
  \Gamma({\cal O}_V( (1+(r+2)+(r+2)^2...+(r+2)^s)D)$$ 
is an incomplete
linear system. 
Let $W_s$ be the rational image
of $V$ in the corresponding projective space. The supremum of
the  dimensions 
of the $W_s$ Is at most the Iitaka dimension of $D.$ 
\m The rational map $V\to W_s$ lifts to a lowest morphism $V_{s+1}\to W_s,$ 
and  the tower of Nash blowups
$...V_2\to V_1 \to V_0=V$ is induced from
the tower of projective (but not birational) morphisms
$...W_2\to W_1 \to W_0\to point.$
\m Stabilization of the tower of $V_s$ is equivalent to
finite type for the graded  sheaf $R =R({\cal O}_V(H))$ while stabilization
of the tower $W_s$ is equivalent to finite type for
the graded algebra of global sections.   It is a triviality (from
both the algebraic and geometric point of view) that stability for
the $W_s$ is equivalent to stability for the $V_s.$ 
\m When $D$ is finitely-generated of Iitaka dimension one,
the graded ring of global sections which controls
finiteness of $...V_2\to V_1\to V_0$ is a sub-algebra of
a finitely-generated algebra of dimension one. Such a sub-algebra
is finitely generated, and so eventually $V_{i+1}=V_i$ and $W_{i+1}=W_i.$ 

\m I should also remark,  the result 
has no reasonable application to the case of Nash blowups, i.e., to the
case when $r=dimension(V),$ because [18] 8.5.5 shows that 
$D-H$ is basepoint free, so $D$ is very ample, and the Iitaka dimension
of $D$ is $r$ itself, for any normal projective variety $V$ which
is Gorenstein with at worst isolated irrational singularities.

\vfill\eject \centerline {2. Definition of $F$} \m
\m In this section we will define the functor $F.$
 Let  
 $V$ be an irreducible variety over a field $k.$ 
Choose a $K$
linear sub Lie algebra (actually any  subspace will do if one is
willing to consider distributions rather than foliations)
${\cal L}\subset
{\cal D}er_k(K,K)$ where $K$ is the function field of $V.$ 
This corresponds
to a singular foliation on $V$. This and other basic definitions can
be found in [1].
 Let $$r=dimension_K({\cal L}).$$ 
\m
For coherent sheaves $A, B$ on $V$  we  define the product
$$AB \buildrel \hbox{def} \over = (A \otimes B)/torsion.$$
\m Let $\Omega$ be  sheaf of differentials
along the foliation, ie  the image of the natural evaluation map to the
dual vector space $\widehat {\cal L}.$
$$ev: \Omega_{V/k}\to \widehat{\cal L}$$
$$ dx \mapsto (\delta \mapsto \delta(x)).$$
for $x$  a local section of the struture sheaf of $V.$
\m For any torsion free rank one coherent sheaf $I$ let $P(I)$
be the sheaf of first principal parts of $I$ with respect to the
foliation. Thus $P(I)$ is  the middle term of the
sequence
$$0 \mapsto I\Omega \to P(I) \to  I \to 0$$
defined by $\alpha=ev_*(\beta)$ where  
$\beta \in Ext^1(I, I\otimes \Omega_{V/k})$
is the Atiyah class of $I.$  Now let
$$F(I) = \Lambda^{r+1}P(I)/torsion.$$
Note $F$ defines a functor which acts on the full subcategory of
torsion free rank one coherent sheaves on $V.$
\m In [1] we defined a less natural operator called ${\cal J}$ acting on
 fractional ideals (ie nonzero finitely-generated submodules of the $R$
module $K$), which depended
 on a generating basis $\delta_1,...,\delta_r\in {\cal D}er_k(K,K)$ 
of the foliation. If $V$ is affine and $I$ is an ideal
in the coordinate ring, the fractional ideal ${\cal J}(I)$  
is the one generated by determinants of matrices whose rows
are $(f, \delta_1f,...,\delta_rf)$ for $f\in I.$ The span of such
rows themselves is a copy of $P(I)$
 and the morphism $P(I)\to I$ 
in the exact sequence above sends this row to its first entry.

\m \vfill\eject\centerline{3. Definition of $R$}
\m Now that we have defined multiplication of coherent sheaves and 
the functor $F$, define torsion free rank one coherent sheaves $J_i$ and $L_i$  for $i=0,1,2,3,...$ by the
inductive rules
$$\left\{\matrix{J_0=I \cr
J_{i+1}=IL_0L_1...Li\cr
L_i=F(J_i)}\right. \ \ \ \ \ (1)$$ These depend functorially on
$I$ and like $L_0=F$ can be viewed as functors acting on 
rank one coherent sheaves.
Also
for $i=0,1,2,,...$ define the rank one torsion free coherent sheaf $R_i$ 
by the rule
$$\left\{\matrix{R_0={\cal O}_V\cr R_i =  L_0^{a_0}L_1^{a_1}...L_s^{a_s}}\right.\ \  \ \ \ \ \ (2)$$
where the numbers $a_\alpha$ are chosen so that the base $(r+2)$
expansion of $i$ is 
$$i=a_0+a_1(r+2)+...+a_s(r+2)^s$$
with $0\le a_\alpha < (r+2).$
\m  {\bf 1. Lemma} If the characteristic of $k$ is zero there is  a multiplication map for each $i,j$  
$$R_i R_j \to R_{i+j}$$
which makes 
$$R=R_0 \oplus R_1 \oplus ...$$
into a 
graded sheaf of ${\cal O}_V$ algebras. 
\vfill\eject\noindent 
\m  Proof of Lemma 1.
  First suppose $V$ is affine, with
affine coordinates $x_0,...,x_n.$ To make our formulas work
  assume $x_0$ is the constant function $x_0=1.$ 
Also choose a $K$ basis
$\delta_1,...,\delta_r\in {\cal L}.$
When $I\subset K$ is the fractional ideal generated by
a sequence of rational functions $y_0,...,y_m$  the isomorphim 
between $P(I)$ and the span of the   rows $(f, \delta_1f,...,\delta_rf)$
induces by passage to highest exterior powers  an embedding
$$F(I)\subset K$$ with image the fractional ideal ${\cal J}(I)$ generated
by determinants 
$$h = determinant\left(\matrix{f_0&\delta_1f_0&...&\delta_rf_0\cr
                               &            &...\cr
                      f_r &\delta_1f_r&...&\delta_rf_r}\right)
\ \ \ \ \ (3)
$$
where $f_0,...,f_r$ run over all pairwise products 
$x_iy_j.$  An observant reader would object that the definition
of ${\cal J}(I)$ on the previous page should require the $f_i$ to
run over all the elements of $I,$ 
however
because of  [1] Proposition 11, and because  $x_0=1$  
 the determinants displayed above generate
the whole of ${\cal J}(I)$ when the $f_i$ run over the smaller
set of pairwise products.

\m 
The following calculation of $h^{r+2}$ can be deduced from  (3)
$$h^{r+2}=determinant
\left(\matrix{hf_0&\delta_1(hf_0)&...&\delta_r(hf_0)\cr
                               &            &...\cr
                      hf_r &\delta_1(hf_r)&...&\delta_r(hf_r)}\right).$$
It is proved as follows: by 
multilinearity of the determinant we can argue as if the $\delta_i$
commute with $h,$ the commutators of the $\delta_i$ and
$h$ cancel out, and so the right side of this
equation   is $h^{r+1}$ times the right side
of (3). The left
side of (3)  is  just $h.$  So the determinant evaluates to
the product $h^{r+1}\cdot h = h^{r+2}.$ This argument is central
to both [1] and the current paper.
\m The `polarization
identity' in linear algebra  then
expresses  any homogeneous  polynomial   of degree $r+2$  over the
rational numbers as a linear combination of such powers $h^{r+2}$ where
$h$ is a linear form. This proves any homogeneous polynomial
of degree $r+2$ in ${\cal J}(I)$ belongs to ${\cal J}(I{\cal J}(I)).$
In this way we obtain an a priori non-natural map  
$$F(I)^{r+2}\cong {\cal J}(I)^{r+2}\to {\cal J}(I{\cal J}(I))
\cong  F(IF(I)).$$ It remains now to verify naturality.

\m Let's now show the definition is independent of choice
of the functions $x_i$ and the rational functions $f_j.$  A 
simpler  but infinite formula not using either choice would have given the
same 
answer as (3) if one   allows the 
$f_i $ on the right side of the equation 
to run over all elements of the fractional ideal  $(y_0,...,y_m)$ and
not only the products $x_iy_j$
(see again [1] Proposition 11). This shows the choice of $x_i$ and $f_j$ is
inessential.  
\vfill\eject\noindent
We now need to show the map is indepenent of choice of basis of ${\cal L}.$ 
If  $\tau_1,...,\tau_r$ were a different basis then
$$\tau_i=\sum_j a_{ij}\delta_j$$ 
with $a_{ij}\in K.$
Then the images of $F(I)^{r+2}$ and $F(IF(I))$ in $K$ are
both  multiplied by  the same rational function
$$det(a_{ij})^{r^2+3r+2}.$$
Because the map is independent of all the choices it then
patches on affine open pieces and defines a map in the case $V$ need not
be affine.
\m Next use the  map $L_i^{r+2} \to L_{i+1}$ to
define a  multiplication  map $R_iR_j \to R_{i+j}.$ For this we 
apply definition
of the $R_i$ in (2).   For example
 if  $r=8$  Then
$(r+2)=10$ and we are using  familiar base ten expansions. 
If we wanted to calculate $98+3$ in the familiar way in base ten, we would
`carry' twice, writing
$$(9\times 10 + 8) + 3 = (9+1)\times 10 + 1 = 10^2+1.$$
Correspondingly we define
$$R_{98}R_3\to R_{101}$$  to be the composite
$$R_{98}R_3=(L_1^9L_0^8)(L_0^3)\to L_1^{10}L_0\to L_2L_0=R_{101}$$
where the first map is indued by $L_0^{10}\to L_1$ and the second by $L_1^{10}\to L_2.$

\m To check the multiplication is well defined and satisfies the ring 
axioms 
it suffices to treat the affine 
case. Then $R$ with its multiplication law has the structure
explicitly described in the next  section. QED
\m \vfill\eject\centerline {4. Elementary Construction of $R$.} 
\m When $V$ is affine with coordinate ring $R_0,$  $\delta_1,...,\delta_r$
are a $K$ basis of ${\cal L}$  and $I\subset K$ is a fractional ideal we
can
define $R$ in an elementary way
to be be the smallest subring  of $K[T]$
  such that
\begin{enumerate}\item[i)] $R$ contains   $R_0\subset K $
in degree zero
\item [ii)]
Whenever the fractional ideal $IR = I\oplus IR_1 \oplus IR_2 \oplus ...$
contains homogeneous 
elements
$f_0,...,f_r$  of the same degree $i$ where  $i$ is either
 zero or equal to a possibly zero partial sum of the divergent geometric series 
$$1+(r+2) + (r+2)^2 +(r+2)^3+...$$
then $R$   must  contain the product
$Td$ where $d$ is the determinant already displayed in (3) 
(now disregarding
the phrase of text following the display).
\end{enumerate}
\m {\bf 2. Remark}  The construction of $R$  appears mysterious
for two  reasons. Not only is the connection with the
base $(r+2) $ expansions unusual, also the choice of degrees for
the grading is strange.Given that $L_0,...,L_i$ have degrees
$1, (r+2), ..., (r+2)^i$  and $F(IL_0L_1...L_i)=L_{i+1}$ has
degree $(r+2)^{i+1}$ we see from the formula
$$(r+2)^{i+1}=1+(r+1)(1+(r+2)+...+(r+2)^i)$$
that the operation  of $F$ has the effect of  multiplying
degrees by $(r+1)$ and then  adding one. So the effect of
$F$ on degrees is neither addition or multiplication, but an
affine transformation. This had mystified me for more than
a  year now, but it now seems the mystery will
be solved in a later section of this paper, when we introduce
separately a very ample divisor $H$ and a canonical divisor $K.$ 
What one shall see is the effect of multiplying an element
in the span of these two divisors by the integer 
$r+1$ and then adding $K.$

\m\vfill\eject \centerline{5. Statement of Results}
\m {\bf 3. Theorem}
 Let $V$ be an irreducible variety over a field $k$ of characteristic zero.
Let ${\cal F}$ be a singular foliation on $V.$ Let $I$ be a sheaf of
ideals on $V.$ Let $V_0\to V$ be blowing up $I$ and  subsequently let
$$ ... V_2 \to V_1 \to V_0 \to V$$
be the Gauss process for the foliation lifted to $V_0.$
The following are equivalent
\begin{enumerate}\item
The sheaf of graded ${\cal O}_V$ algebras $R$ defined using
${\cal F}$ and $I$ is finite type over $k$
\item The tower $...V_2 \to V_1 \to V_0 \to V$ stabilizes
\item The natural map $F(J)^{r+2}\to F(JF(J))$ is an isomorphism
for some $J=IL_0...L_{t-1}$ and some $t\ge 0.$
\end{enumerate}
\m We will give two separate proofs of the theorem later on.

\m {\bf  4. Corollary.} A foliation on an irreducible variety $V$ in characterist
ic zero
can be resolved by a locally projective birational  map (i.e., a blow up of a sh
eaf
of ideals) if and only if there is an ideal sheaf $I$ on $V$ with
$F(IF(I))=F(I)^{r+2}.$ \m Proof of Corollary. Note that the
 words ``ideal sheaf'' can be replaced by ``torsion free coherent
sheaf of rank one'' as any such sheaf can be  embedded after
twisting by a divisor.
Assume there is an ideal sheaf $I$ with $F(IF(I))=F(I)^{r+2}.$
Then [1] Theorem 15 part ii) with $N=0$ shows that
the Gauss process starting with blowing up $I$ finishes
in one further step.
Conversely
assume the foliation can be resolved
and let   $V_0$ be a resolution, with
resolving ideal sheaf $I.$ This time the
Gauss process for the foliated variety $V_0$ finishes at step 0.
Though there is
no direct connection between the number of steps in the Gauss resolution (zero i
n this case)
 and the
generating degrees of $R,$ nevertheless, by part 1. of Theorem 3 we do
know $R$ is finite type.
Then by part 3. of Theorem 3 the sheaf $J=IL_0...L_{t-1}$ satisfies  $F(JF(J))=F
(J)^{r+2}$ for
some  finite value of $t.$ QED
\m As we explined in the introduction, one application
is this
\m {\bf 5. Corollary} Let $V$ be 
a normal foliated, irreducible projective
variety  
over an algebraically closed field $k$ of characteristic zero
such that $(r+1)H+K$ has Iitaka dimension  less than two,
where $H$ is a hyperplane section of $V,$ $r$ is
the dimension of the foliation and $K$  is the canonical
divisor of the Foliation.
Then the sequence of  Gaussian  blowups of $V$  stabilizes.
\m {\bf Remark.} Because the description of generating
sections of $L_0 \subset {\cal O}_V(K+(r+1)H)$ always 
includes a nonzero global section,
 There is always a choice of  effective divisor
linearly equivalent to $(r+1)H+K,$  and the Iitaka
dimension is  the transcendency degree of the algebra
of global sections of the structure sheaf of the quasi
projective variety $V\setminus D.$

\m
\vfill \eject  \centerline{6. Four affine  examples.}
\m{\bf 1.  Example.} Firstly consider the cusp, whose
coordinate ring is $k[x^2,x^3, x^4,...].$ The one dimensional
 foliation is spanned by any nonzero derivation so
let's use  $x\partial/\partial x.$
Take $I$ to be the unit ideal.  Our ring $R$
is now the smallest subring of $k[x,T]$ containing $k[x^2,x^3,...]$
and with the property
that whenever $R$ contains two monomials $A,B$ whose $x$ degrees
are distinct and whose $T$ degrees are the same, possibly
zero,  partial sum
of the divergent 
  geometric series
$1+3+3^2...$ then the product $ABT$ must be contained in $R.$ 
Here is therefore the list of monomials of low degree in $R.$
\m $$\matrix{x^8&Tx^8&T^2x^8&T^3x^8&T^4x^8\cr
x^7&Tx^7&T^2x^7&T^3x^7&T^4x^7\cr
x^6&Tx^6&T^2x^6&T^3x^6\cr
x^6&Tx^5&T^2x^5&T^3x^5\cr
x^4&Tx^4&T^2x^4\cr
x^3&Tx^3
\cr x^2&Tx^2\cr\cr 1}$$
\m The monomials $Tx^2$ and  $T^3x^5$ are the only interesting ones. They
are the only nontrivial monomials in $R$ which are not a product of monomials of
smaller $T$ degree.
$T^3x^5$ is included  because $1$ is a partial sum of a geometric series of the
base 3,  and
the separate monomials $Tx^2$ and $Tx^3$  have distinct $x$ degrees,
and so their product times $T$ is included.
 The monomials in $L_iT^{3^i}$ are the columns in
the diagram which are indexed by powers of three, and the last observation
implies the ninth column
is the third power of the third column.
Then  $L_1^3=L_2$
so the desired equation $F(J)^3=F(JF(J))$  holds
for $J=IF(I).$  This example and the next one are included as 
contrasting tutorial examples. 
\vfill\eject\noindent {\bf 2.  Example.} For a second example consider the 
singular  foliation
of the plane given by the vector field 
$$x \partial/\partial x + 2 y \partial/\partial y.$$
It is known that this foliation cannot be 
resolved by blowing up points, and that singular foliations
of the plane can rarely be resolved by blowing up points. 
There is a 
general `desingularization theorem' for one dimensional 
singular foliations of the plane, which for example is applied to
complex singularities which are defined over the reals, as described
in Ilyashenko's
centennial history of Hilbert's 16'th problem [17]. According
to Ilyashenko's article, the 
desingularization theorem has a long history,
with contributions by Bendixson, Seidenberg, Lefshetz, Dumortier,
and van Essen, but is not written up in any book; it states that
after blowing up points, a one dimensional foliation of the plane 
can be arranged to have isolated singular points with nonzero
linear part, such that  the linear part at each singular point has
a nonzero eigenvalue. The question of which singular foliations
can actually be resolved in the stronger sense of lifting to a
nonsingular foliation after pulling back (and taking
an irreducible component) by a locally projective
birational morphism, even for linear foliations, is unsolved
as far as I know and will not be solved in this paper, though 
we will consider which foliations can be resolved torically, and
we expect there to be no surprises if resolutions which are
not toric are allowed. If we apply our results to the case of
linear plane foliations we will see that the  linear
foliations of the plane which can be resolved torically
with respect to an eigenspace  decomposition have the
eigenvalue pairs
$(1,0),(0,1)$ and $(1,1).$  

\m 
For the linear foliation of the plane we are considering now, with
the eigenvalues $1$ and $2,$ it is an easy calculation that
the foliation cannot be resolved by blowing up points, and by 
Zariski  [4]
any 
proper birational map of normal surfaces arises that way. 
It is therefore not possible 
that this particular foliation
could be resolved by a proper birational map from a normal surface,
and this is easily seen directly.
In this  example and the next one, we will try to understand
the deeper reason why it cannot be resolved, to motivate
a more general theorem.  In this example
let us look at why the Gauss process itself does not converge.
In view of the main theorem, it shall be the same thing to
 examine the ring $R$ and verify that  it is not finitely generated.
 A $k$ basis of $R$ consists of the
 the smallest set of  monomials in $x,y,T$ which
is closed under multiplication,  such that
all monomials $x^iy^j$ for $i,j\ge 0$ are included, and in addition
if $A$, $B$ are two   monomials  in $x$ and $y$ of distinct degree
(where $x$ is given degree 1 and $y$ degree 2),  and if $AT^i$ and $BT^i$
are included where $i$ is a partial sum of the
geometric series $1+3+3^2 +...$
then so is $T\cdot AT^i \cdot BT^i=T^{2i+1}AB.$
 Consider the monomials of the form  $x^vy^jT^k$  which
occur.
For $v=0$ we get  $y^jT^k$ with $j\ge k \ge 0.$
These arise as powers of $y$ times powers of $yT.$
For $v=1$
we get the smallest  set of monomials containing
$x$ and closed under multiplication by both $y$ and $yT$
and  under 
$$xy^jT^i \mapsto xy^{j+i}T^{2i+1}     $$
when $i$ is zero or of the form $1+3+3^2...$.
The latter rule comes from the product $T\cdot (m_0T)\cdot (m_1)T$
where $m_0=y^i$ and $m_1=xy^j,$
as note that $xy^j$ always has odd degree which
is always distinct from that of
$y^i$, which is even.
Let
 $P(x,y,T)$ be the sum of all monomials  in $R.$ Then
${\partial \over {\partial x}} P(x,y,t) |_{x=0}$ is the sum
of the smallest set of monomials in $y $ and $T$ which
contains  all $y^{3^i-i-1}T^{3^i}$  for $ i=0,1,2,...$
and  is closed under
multiplication by $y$ and $yT$. This is
$${1\over{1-yT}}(
1+{{yT}\over y}+{{(yT)^3}\over{y^2}} + {{(yT)^9}
\over {y^3}}+{{(yT)^{27}}\over {y^4}}+...).$$
which is not a rational function so  $R$ is not
finitely generated. 
\vfill\eject\noindent {\bf 3.  Example.} For the third example, let us continue
to consider the
same foliation as in 
the previous example, but now let us 
try to understand why there is no toric resolution of this one-dimensional
 foliation.
Our proof will introduce a technique which will be able to generalize
to foliations of higher dimension.
A toric resolution is a special case of 
of  blowing up a monomial ideal $I.$ We shall show in fact that
the foliation is not resolved by blowing up any monomial ideal, and
in fact (which is equivalent) that it cannot be resolved by blowing
up a monomial ideal and then following by a finite chain of Gauss blowups.
Consider the infinite sequence of maps 
$$... V_2 \to V_1 \to V_0\to V$$ 
 where $V_0\to V$ is blowing up along $I$ 
and $V_{i+1}\to V_i$ is the Gauss blowup
along the foliation for $i\ge 1.$  
I claim the process can not converge, meaning none of the maps
is an isomorphism. The theorem again says this is 
equivalent to finite generation of the appropriate ring $R.$ 
Like in the previous example, finite generation of $R$ would correspond
to an equation $L_i^3=L_{i+1}.$ Writing $J=IL_0...L_{i-1}$ we
would have $F(J)^3=F(JF(J)).$  Let $A$ be the smallest monomial
in $J$ for the alphabetic ordering (where $x$ is given more priority
than $y$). Let $B=Ay$ and $C=Ay^{-e}x$ with $e$ maximum. The
the monomials $AB$ and $AC$ belong to $F(J)$  and so $A^2B$ and
$A^2C$ belong to $JF(J).$ The degrees have opposite parity so
are unequal, whence $(A^2B)(A^2C)$ belongs to $F(JF(J)).$ 
The product cannot be rewritten any other way and it cannot
belong to $F(J)^3$ because one of the three factors would have to be
$A^2$ which does not belong to $F(J)$ at all.
\vfill\eject\noindent {\bf 4.  Example.} Our fourth example
involves examining  the conditions for resolvability of the 
unique codimension zero foliation on a normal  irreducible affine  complex 
algebraic 
 surface.  
Taking affine coordinates $x_0,...,x_n,$ we assume  $x_0$ is the
constant function $x_0=1$ and that $x_1,x_2$ are algebraically
independent. Let $y_0,...,y_m$ be
rational functions and let $V_0=Bl_IV$ where $I$ is the fractional ideal
generated by $(y_0,...,y_m).$  
\m  We
choose as our derivations $\partial/\partial x_1$ and
$\partial/\partial x_2.$ 
 Since the characteristic
of $k$ is zero, and any rational function on $V$ is algebraic
over $k(x_1,x_2)$,  and so the $\partial/\partial x_i$
for $i=1,2$ can be evaluated on any rational function, and thus
define a pair of (commuting) derivations on the function field. 
For any fractional ideal $I$ we may use these two derivations to
view $F(I)$ again as a fractional ideal, as explained in 
section 4. The map $F(I)^{r+2}\to F(IF(I))$ is then an
an inclusion of fractional ideals, depending on the rational functions
$y_0,...,y_m.$  
\m Because the generating  sequence of $I$
is a product with the sequence $x_0,...,x_n$
with $x_0=1$ the hypothesis of [1] Proposition 11  is satisfied and
shows $F(I)$ is generated by
determinants

$$
  \left|\matrix{
x_ey_f&\partial/\partial x_1(x_ey_f)&\partial/\partial x_2(x_ey_f)\cr
x_gy_h&\partial/\partial x_1(x_gy_h)&\partial/\partial x_2(x_gy_h)
}\right|
$$
 Since the generators of $I$ satisfy the hypothesis of the
proposition so do these generators times the determinants displayed
above. Applying the same propositionto $IF(I)$ shows 
$F(IF(I))$ is generated by the determinants

\m

\noindent {\resizebox{6.0 in}{1.0 in}{
$ \left| \matrix{
             x_ay_b\left|\matrix{
x_cx_d&\partial/\partial x_1(x_cy_d)&\partial/\partial x_2(x_cy_d)\cr
x_ey_f&\partial/\partial x_1(x_ey_f)&\partial/\partial x_2(x_ey_f)\cr
x_gy_h&\partial/\partial x_1(x_gy_h)&\partial/\partial x_2(x_gy_h)
}\right|
&
\partial/\partial x_1
\left(
x_ay_b\left|\matrix{x_cy_d&\partial/\partial x_1(x_cy_d)&\partial/\partial
x_2(x_cy_d)\cr
x_ey_f&\partial/\partial x_1(x_ey_f)&\partial/\partial x_2(x_ey_f)\cr
x_gy_h&\partial/\partial x_1(x_gy_h)&\partial/\partial x_2(x_gy_h)
}\right|     \right)
&\partial/\partial x_2
\left(
x_ay_b\left|\matrix{x_cx_d&\partial/\partial x_1(x_cy_d)&\partial/\partial
x_2(x_c
y_d)\cr
x_ey_f&\partial/\partial x_1(x_ey_f)&\partial/\partial x_2(x_ey_f)\cr
x_gy_h&\partial/\partial x_1(x_gy_h)&\partial/\partial x_2(x_gy_h)
}\right|     \right)
\cr \hbox{\phantom{THIS}}\cr              x_iy_j\left|\matrix{
x_kx_l&\partial/\partial x_1(x_ky_l)&\partial/\partial x_2(x_ky_l)\cr
x_oy_p&\partial/\partial x_1(x_oy_p)&\partial/\partial x_2(x_oy_p)\cr
x_qy_s&\partial/\partial x_1(x_qy_s)&\partial/\partial x_2(x_qy_s)
}\right|
&
\partial/\partial x_1
\left(
x_iy_j\left|\matrix{x_kx_l&\partial/\partial x_1(x_ky_l)&\partial/\partial
x_2(x_k
cy_l)\cr
x_oy_p&\partial/\partial x_1(x_oy_p)&\partial/\partial x_2(x_oy_p)\cr
x_qy_s&\partial/\partial x_1(x_qy_s)&\partial/\partial x_2(x_qy_s)
}\right|     \right)
&\partial/\partial x_2
\left(
x_iy_j\left|\matrix{x_kx_l&\partial/\partial x_1(x_ky_l)&\partial/\partial
x_2(x_k
y_l)\cr
x_oy_p\beta&\partial/\partial x_1(x_oy_p)&\partial/\partial
x_2(x_oy_p)\cr
x_qy_s\delta&\partial/\partial x_1(x_qy_s)&\partial/\partial
x_2(x_qy_s)
}\right|     \right)
\cr \hbox{\phantom{THIS}}\cr
             x_ty_u\left|\matrix{
x_cx_d&\partial/\partial x_1(x_cyd)&\partial/\partial x_2(x_cy_d)\cr
x_ey_f&\partial/\partial x_1(x_ey_f)&\partial/\partial x_2(x_ey_f)\cr
x_gy_h&\partial/\partial x_1(x_gy_h)&\partial/\partial x_2(x_gy_h)
}\right|
&
\partial/\partial x_1
\left(
x_ty_u\left|\matrix{x_cx_d&\partial/\partial x_1(x_cyd)&\partial/\partial
x_2(x_cy_d)\cr
x_ey_f&\partial/\partial x_1(x_ey_f)&\partial/\partial x_2(x_ey_f)\cr
x_gy_h&\partial/\partial x_1(x_gy_h)&\partial/\partial x_2(x_gy_h)
}\right|     \right)
&\partial/\partial x_2
\left(
x_ty_u\left|\matrix{x_cx_d&\partial/\partial x_1(x_cyd)&\partial/\partial
x_2(x_c
y_d)\cr
x_ey_f&\partial/\partial x_1(x_ey_f)&\partial/\partial x_2(x_ey_f)\cr
x_gy_h&\partial/\partial x_1(x_gy_h)&\partial/\partial x_2(x_gy_h)
}\right|     \right)
} \right|  $
}}
\m 
\m 
The equation $F(I)^4=F(IF(I))$ holds just when the large determinants
can be expressed as  
homogeneous polynomials of degree four in the small deterimants.
Because $V$ is normal, the coefficients of the
degree four polynomials whenever they exist
would restrict to well-defined
 functions on the regular locus of $V.$ 
Likewise,
all $y_i$ in the equation can be  multiplied
by a common denominator 
$t$ of the coordinate ring,  
both sides of the equation are
multiplied by $t^{(r+1)(r+2)}.$ Then too the $y_i$ may be considered
as well-defined functions on the whole of the regular locus of $V.$
The equation is therefore equivalent to
a system of partial differential equations on a complex
manifold. This is the system of differential equations one  would
have encountered in the real case if one were to have  looked 
for an algebraic condition
which bounds the Gaussian curvature of the leaves under a projective
embedding.
A solution $(y_0,...,y_m)$ need not generate
a resolving ideal, nor does a set of generators of a resolving
ideal necessarily provide a solution. 
Zariski's work [3] showing surfaces can be resolved
does  provide a resolving ideal. The chain  of subsequent
Nash blowups finishes because it 
is trivial. The ring $R$ is an algebraic
model of the chain, which however need not become trivial
at the first step.  Theorem 3 part 1
 shows $R$ is finite type.
Then there is an ideal $J$ satisfying part 3 of the same theorem. 
Taking $y_0,...,y_m$ to be a generating sequence we see therefore that
this, possibly familiar, system of differential
equations  on the smooth manifold of $V$  always has a solution.
\vfill\eject\centerline{7. Foliations coming from split semisimple representations of
commutative  Lie algebras.}
\m The simplest examples of singular foliations arise from a faithful
split semisimple representations of an 
$r$ dimensional   commutative  Lie algebra $G$ on a
vector space $V$ (all over the base field $k$ of characteristic zero). We 
choose coordinates $x_1,...,x_n$ on $V$ and elements in the dual
of $G$
$\alpha_1,...,\alpha_n\ \in \widehat G$
so that 
 for $s \in G$ we  have
$$sx_i=\alpha_i(s)x_i.$$
We may view the vector space 
$V$ as a toric variety such that the monomials in $x_i$ 
are the characters of the torus which are well-defined on all of $V.$
\m {\bf  6. Corollary.}  It is possible to resolve the foliation
on   $V$ torically if and  only if the underlying set of nonzero roots
(counted without multiplicities) forms a basis of the dual Lie algebra.
\m  {\bf  7.  Remarks.}
The proof shows more. It allows a more general blowup of
an arbitrary monomial ideal.
  Also
the theorem  applies in cases when the blowup may
be a singular variety.
\m Proof. For any monomial ideal $I$ 
view $F(I)$ as the ideal  generated
by  determinants of matrices whose rows
are the $(f, \delta_1 f , ... \delta_r f)$
where $f,g\in I.$ By multilinearity of the determinant,
 we may restrict $f$  to run over monomials in $I$, and so
again $F(I)$ is a monomial ideal. In this way we inductively
see that all $L_i$ are monomial ideals.
\m Suppose first that the foliation can be resolved by some
monomial ideal $I$ then Theorem 3 implies 
 that    $F(J)^{r+2}=F(JF(J))$ for some $J=IL_0...L_{t-1},$
again a monomial ideal.
Let $A$ be the monomial in $J$  which is minimum for the alphabetical
ordering on $(i_1,...,i_n)$ where $i_1 $ is given least significance. Also
let $B$ be the smallest monomial in $J$ in which the power of $x_{r+1}$
is one larger than what occurs in $A.$
\vfill\eject\noindent Since the $V$ is a faithful representation of
our Lie algebra $G,$  
we know there must exist a basis of   $\widehat G$ consisting of roots.
By  choice of numbering these may be assumed to be
 $\alpha_1,...,\alpha_r.$  If $n=r$ we are done, so assume
$n\ge r+1.$
\m For any monomial 
$M=x_1^{i_1}...x_n^{i_n}$  write $f(M)=i_1\alpha_1+...  +i_n\alpha_n.$
In this way we obtain a homomorphism from the set of monomials
(characters of the torus that are well defined on all of $V$)
to the dual of $G.$ 
Consider the product $A^{(r^2+3r+1)}(x_1...x_{r-1})^{r+2}x_r^{r}B.$
First we associate this as  a product of $r+1$ factors
$$(\prod_{i=1}^{r-1} [x_iA^{r+2}(x_1...x_r)])\cdot [A^{r+2}(x_1...x_r)]
\cdot [ A^{r+1}x_1...x_{r-1}B]$$
Let us show 
that each  factor is  a product of $r+2$ monomials with $f$ values which 
span $\widehat G$ affinely.
The first factor is
$$A \cdot A \cdot Ax_1^2 \cdot Ax_2 \cdot ... \cdot A x_r$$
and applying $f$ to each term and subtracting $f(A)$ yields
$$0,0, 2\alpha_1, \alpha_2, \alpha_3, ..., \alpha_r$$
which  affinely span since the $\alpha_1,...,\alpha_r$ are
linearly independent. All but the last factor behave
 similarly to this one with the coefficient of $2$
occuring in a different position  or being absent.  The last factor
is $A\cdot A \cdot   Ax_1 \cdot ... \cdot Ax_{r-1} \cdot B$ This time
after subtracting $f(A)$ from each term, the sequence of $f$ values
is $0,0,\alpha_1,...,\alpha_{r-1},f(B)-f(A).$
Modulo the hyperplane spanned by the other terms, 
the last term is congruent to 
to $\alpha_{r+1}-i\alpha_r$ for some positive integer $i.$ We shall
aim to prove  this is congrent to zero for some positive value of $i.$
 Since this deduction will be true after any permutation of the
variables it will be a strong restriction, which will imply our conclusion.
In the way of hoping for a contradiction, then, suppose that
there is no positive number $i$ so that  $\alpha_{r+1}$ is
 congrent to $i\alpha_r$  modulo the span of $\alpha_1,...,\alpha_{r-1}.$ 
Then the final factor too 
is a product of $(r+1)$ monomials in $J$ whose images in $\widehat G$
are affinely independent.
\m Next let us  apply $f$ to each of the $(r+1)$  whole factors in
square brackets.
Subtractiing $(r+2)f(A)+\alpha_1 + ... + \alpha_r$ from each we obtain
the sequence
$$\alpha_1, \alpha_2, ..., \alpha_{r-1},0, f(B)-f(A)-\alpha_r.$$
The last term is congrent to $\alpha_{r+1}-(i+1)\alpha_r$ and 
again using our assumption that this is not in the span
 of $\alpha_1,...,\alpha_{r-1}$ we have a sequence
that is affinely independent. 
This establishes that our monomial belongs to
$F(JF(J)).$
\m Because of our assumption that $F(JF(J))=F(J)^{r+2}$  
 our monomial must then belong to $F(J)^{r+2}$. This means  
it must be possible to refactorize our product  a different way, as a product
 of $(r+2)$  monomials in $J$ so that each monomial
factorizes as a product of $r+1$ monomials with affinely independent
images in $\widehat G.$ 
The most significant letter where our product differs from
$A^{(r+2)(r+1)}$ is in the exponent of $x_{r+1}$, which is  one larger.
In our factorization as a product of $(r+2)(r+1)$ monomials in $J,$ 
choose one of the monomials which includes a higher power of $x_{r+1}$
than $A$ does. Remove  this  factor from our product expression, 
along with whichever $r$ further factors  
are associated with it.
At most  $r$ of the remaining $(r+1)^2$
factors can have a power of $x_r$ that is larger than what  occurs
in $A.$ Applying $f$ to each factor we 
obtain a sequence of  $(r+1)^2$ elements of $\widehat G$ of which  at least
$r^2+r+1$ lie in an affine hyperplane, let us call it $H.$
The sequence is a disjunction of 
$r+1$  subsequences,  each of which consists of 
$(r+1)$ affinely independent elements.  An affinely independent set
always 
contains an element in $\widehat G \setminus H$  and so there are at 
least $r+1$ terms of our sequence in $\widehat G \setminus H.$ 
This is the desired contradiction.
\m We have shown, therefore, that for every set of roots $\alpha_1,...,\alpha_{r+1}$
such that $\alpha_1,...,\alpha_{r-1}$ are linearly independent, it must be the
case that  $\alpha_r-i\alpha_{r+1}$ lies in the linear 
span of $\alpha_1,...,\alpha_{r-1}$ for some integer $i\ge 0. $  
Interchanging $r$ and $r+1$ we see that this is true for $i=1.$ Letting
$\alpha_r$ and $\alpha_{r+1}$ range over all roots not equal to $\alpha_1,...,\alpha_{r-1}$
we see that all roots belong to the hyperplane spanned by $\alpha_1,...,\alpha_{r-1}$
together with one affine translate of that hyperplane. Applying this principle
to a basis $\alpha_1,...,\alpha_r$ of $\widehat G$ we see that all roots 
are of the form $a_1\alpha_1 + ... + a_r\alpha_r$ for $a_i \in \{0,1\}.$ 
So any two sets of  roots which are  vector space bases of  $\widehat G$ 
give rise to a change-of-basis matrix with positive integer entries. But the
only invertible elements in the monoid of matrices 
with positive integer entries  are  permutation
matrices. QED

\m \vfill\eject\centerline { 8. Proof of Theorem 3.}
\m In this paper we will give two proofs of the theorem.
The first proof in this section does not use any algebraic geometry.
A subsequent
proof in section 11  will be outlined which is simpler but
uses concepts and theorems of algebraic geometry. 

\m 
 Let us first treat the affine case. Let 
  $V$ be an irreducible affine variety over a field $k$ of characteristic
zero  with
function field $K,$
${\cal L}\subset {\cal D}er_k(K,K)$  a $K$ linear Lie algebra, with
fixed basis $\delta_1,...,\delta_r$.
We will call the coordinate ring of $V$ $R_0$ in the expectation that
it will later become the zero degree component of a graded ring $R.$
Most simply, the torsion free $R_0$ modules
 $L_i, J_i, R_i \subset K$ can be identified with the fractional
ideals ($=R_0$ modules)  
defined by (1),(2), and (3) where in (3) one 
may take $f_i$ to run over all
elements of $I.$
\m
 Step 1: Let us prove 1. $\Rightarrow$ 3 $\Rightarrow$ 2.  in Theorem 3. 
Thus suppose 1. We suppose
 $R$ is finite type. By [1] Lemma 4 and induction,  
the resolution process $...V_2\to V_1 \to V_0 \to V$
has the property that $V_i=Bl_{J_i}(V)$ where 
we recall $J_i=IL_0...L_{i-1}.$ 
By construction,
  $R$ is generated by $R_0$ and  the $L_i$ in degree $(r+2)^i$ 
for $i=0,1,2,....$  By very  elementary properities
of integer base $(r+2) $ expansions, the only way  $R$ can be finite type 
is if actually
$L_{i+1}=L_i^{r+2}$ for some $i.$  This proves part 3. 
Now  by [1] Lemma 4 and induction,
the resolution process $...V_2\to V_1 \to V_0 \to V$
has the property that $V_i=Bl_{J_i}(V)$ where
we recall $J_i=IL_0...L_{i-1}.$
 This implies blowing up $J_{i+1}$ or
$J_{i+2}$ has the same effect so $V_{i+2}=V_{i+1}$ and part 2. is proved.
\m The remaining steps  will be concerned with proving 2. $\Rightarrow$ 1.
in Theorem 3.   From now on we will assume the resolution process finishes.
\m Step 2: In this step we will define a subring of $R.$
We are assuming that 
resolution process finishes in finitely many steps. We know by
[1] Theorem 15 (ii) 
for some $N$ 
$$J_i^NL_i^{r+2}\to J_i^NL_{i+1}\ \ \ \ \ \ \ (4)$$
is an isomorphism.  
\m {\bf Remark.} We alternatively could  now deduce this  
for a possibly smaller value of $i,$ using the fact that
 eventually the Atiyah class on $V_i$ becomes ``multiplicatively trivial''
in the sense that the tensor product of the highest
exterior powers of the  end terms of the 
Atiyah sequence maps isomorphically to the highest
exterior power of the middle term. 
The fact that we know the sheaves
become the same after multiplying by a power of $J_i$ is related to
the fact that the Atiyah classes are calculated on the $V_i$ and not $V.$ 
The rough idea can be summarized as saying now that finite
generation questions on $V_i$ and on $V$ should be equivalent because
 coherence is preserved under direct images of a proper map.

\m   Fix the values of $i$ and $N$ for the
rest of the proof.
\m  Let the graded
subring
$R^{[i+1]}\subset R$ be the subring which is obtained by multiplying
each degree $j$ component $R_j$ by $(IL_0L_1...L_i)^j=J_{i+1}^j.$ We have
$$R^{[i+1]}=\oplus_jJ_{i+1}^jR_j.$$ Before we  prove $R$ is
finite type over $k$  we will  first  
 prove the subring $R^{[i+1]}$ is finite type over $k.$
Our proof will implicitly determine bounds for the degrees
of the generators of $R^{[i+1]}.$  
\noindent
\vfill\eject\noindent Step 3:
In this step we will prove a lemma which will later be useful
in proving that the subring
$R^{[i+1]}$ is finite type.
\m {\bf 8. Lemma} For any $s\ge 0$ 
$$J_{i+1}^{N(r+1)^s}L_{i+s+1}=J_{i+1}^{N(r+1)^s}L_i^{(r+2)^{s+1}}.$$
\m Proof. If $s=0$ multiply the previous formula (4)  by 
$L_i^N$ and use $J_{i+1}=L_iJ_i.$ If $s\ge 1$ assume the lemma
true for smaller $s$ and we have
that the left hand side equals
$$J_{i+1}^{[N(r+1)^{s-1}-1](r+1)+(r+1)}L_{i+s+1}$$
$$=J_{i+1}^{[N(r+1)^{s-1}-1](r+1)+(r+1)}F(IL_0L_1...L_{i+s})$$
Using Theorem 12 of [15] which implies
$(J^{r+1})F(X)\subset F(JX)$ for any $X$ we find
this is 
$$\subset J_{i+1}^{r+1}F(J_{i+1}^{N(r+1)^{s-1}-1}L_0...L_{i+s})$$
Recall $J_{i+1}=IL_0...L_i$ so we know $J_{i+1}^{N(r+1)^{s-1}-1}IL_0...L_i
=J_{i+1}^{N(r+1)^{s-1}}.$
The inductive
hypothesis  implies $J_{i+1}^{N(r+1)^{s-1}}L_{i+\alpha}
=J_{i+1}^{N(r+1)^{s-1}}L_i^{(r+2)^\alpha}$ for $\alpha=0,1,...,s.$    
Simplifying the expression
above using both rules yields
$$=J_{i+1}^{r+1}F((IL_0...L_i)^{N(r+1)^{s-1}}L_i^{(r+2)+(r+2)^2 + ... + (r+2)^s}).$$
This is of the form $J_i^{r+1}L_i^{r+1}F(J_i^aL_i^b)$ for some 
numbers $a$ and $b$ and here the way we proceed depends on whether
$a$ and $b$ are odd or even. If they are both even we use
the fact that $F(X^2)=X^{r+1}F(X)$ for any $X$, which is a special
case of [1] Theorem  14, while if  for 
example $a$ is odd and $b$ is even
 we first apply the inclusion  $\subset L_i^{r+1}F(J_i^{a+1}L_i^b)$
and then apply the rule $F(X^2)=X^{r+1}F(X)$.  In
this way we obtain
$$\subset J_{i}^{N(r+1)^s}L_i^{(r+1)+[N(r+1)^{s-1}+(r+2)+...+(r+2)^s-1](r+1)
}F(J_iL_i).$$
But $F(J_iL_i)=F(J_{i+1}) = L_{i+1}$ can be replaced by $L_i^{r+2}$
because it occurs multiplied by at least $J_{i+1}^N$ and we obtain
$$=J_{i+1}^{N(r+1)^s}L_i^{(r+2)^{s+1}}$$
as claimed. 
\vfill\eject\noindent Step 4: In this step let us just observe
that it follows from lemma 8 by  letting $t=s+i+1$ that
for $t$ suficiently large
$$J_{i+1}^{(r+2)^t}L_t = J_{i+1}^{(r+2)^t}L_i^{(r+2)^{t-i}}
$$ $$ = (J_{i+1}^{(r+2)^i}L_i)^{(r+2)^{t-i}}$$
This is because then  the inequalities
$$({{r+2}\over{r+1}})^t\ge N.$$
and $$s\ge 0$$
will both hold.
\m Step 5. Now we can show that the
subring $R^{[i+1]}\subset  R$ defined at the beginning of the proof
is finite type. 
If  the base $r+2$ expansion of $t$ is 
$$a_0+a_1(r+2) +...+a_m(r+2)^m$$ then the 
degree $t$ component of  of $R^{[i+1]}$ is $$R_t^{[i+1]}
=(J_{i+1}L_0)^{a_0}...(J_{i+1}^{r+2}L_1)^{a_1}(J_{i+1}^{(r+2)^2}L_2)^{a_2}...
(J_{i+1}^{(r+2)^m}L_m)^{a_m}$$ 
and so  step 4 shows $R^{[i+1]}$ is generated by
the 
$J_{i+1}^{(r+1)^t}L_t$ for $$t< max[i+1, {log(N)\over{(log(r+2)-log(r+1))}}].$$ 

\m Step 6.  Now we
can deduce  that $R$ itself is finite type. This
step does not give explicit bounds as it is based on Hilbert's basis theorem.
To simplify notation, let us rename the graded ring 
$R^{[i+1]}$ by the name $W$.
Since
$W$ is finite type over $k$ and therefore 
over $R_0$ it has a sequence of homogeneous
generators $x_1,...,x_t$ of positive degree.  Let $d$ be the
least common multiple of the $degree(x_i)$ and consider all possible
elements of $W$ of degree $d$ which occur as monomials in
the $x_i$. These monomials generate the `truncated' ring
$\oplus_j W_{dj}$.  
The fact that all generators belong to $W_d$ implies
$$W_{dj}=W_d^{j}$$
Using the definition of $W=R^{[i+1]} $ this tells us
that the inclusion
$$R_d^j\subset R_{dj}$$ becomes
an equality after multiplying by
$(IL_0...L_i)^{dj}.$ 
Take $x\in R_{dj}$ and let  $e_1,...,e_n$ be a set of
$R_0$ module generators of $(IL_0...L_i)^{dj}.$ 
We have $xe_i=\sum a_{ij}e_j$ with $a_{ij}\in R_d^j.$ We now apply
Zariski and Samuel's trick of taking deteriminant of the matrix 
$(a_{ij}-x\delta_{ij}).$ 
This gives 
$0 = x^nx^{n-1}b_{n-1}+...+b_0$ with  the $b_t \in R_d^{j(n-t)}.$
We  view this as an equation of integral dependence where
the coefficients $b_t$ lie in the graded ring $\oplus_i R_d^i$
and this shows every element of $R$ whose degree
is a multiple of $d$  is integral over this ring. Since
$R_d$ is finitely generated as $R_0$ module and $R_0$ is finite
type, the ring  $\oplus_i R_d^i$ is finite type over $k$ and by the theorem
of `finiteness of normalization'  we know the truncation
$\oplus_i R_{id}$ since it consists 
only of integral elements, is also finite over
over $\oplus R_d^j.$ Therefore the truncation is finite type over $k$ and
therefore Noetherian.
Finally, the ring $R$ is isomorphic to a
direct sum of $d$  ideals  over the truncated
ring (multiplying by any element of
degree congrent to $-i$ mod $d$  gives an embedding of the
sum of terms congruent to $i$ mod d into the truncated ring).  
Since ideals in  a Noetherian ring are finitely generated, we conclude
$R$ is finite over the truncated ring. It follows then that $R$ is
finite type over $k.$ 
\m Step 7. The proof in the affine case is done, but 
we have to consider the case  $V$ is not affine. Then more functorial definition
of $R$ shows it patches on a suitable affine open cover. QED
\vfill \eject \centerline{9. Discussion}
\m Although the proof of Theorem 3 was  algebraic, it might be helpful
to visualize the corresponding diagram involving the
prescheme $X=Proj(R)$, which does exist even when $R$ is not finite type.
Letting $V_0=Bl_I(V)$ and letting $...V_2 \to V_1 \to V_0$ be the Gauss
process for this variety, then writing 
$X_i=Proj(R^{[i]})$ we  have a pullback diagram
$$\matrix{...&\to&X_2&\to&X_1&\to&X_0&\to&X\cr
            &&\downarrow&&\downarrow&&\downarrow&&\downarrow\cr
          ...&\to&V_2&\to&V_1&\to&V_0&\to&V}.$$
If the Gauss process for $V_0=Bl_I(V)$ converges then
for some $n$ we have $V_{n+1}=V_n$ and $X_{n}=X_{n-1}.$ Then the vertical map $X_{n}\to V_n$
is an isomorphism showing $V_n$ is the pullback
of $V_{n-1}\to V \leftarrow X.$ Thus $X$ is dominated by a variety.

\m\centerline{10. Algebraic geometry interpretation}
\m Let $V$ is a normal irreducible quasiprojective
variety over a field $k$ of characteristic zero, and let $H$ be 
a very ample divisor on $V.$  Choose a
singular foliation  of dimension $r.$ We shall consider the
sheaf ${\cal O}_V(H)$ and the graded sheaf ${ R}({\cal O}_V(H)).$
\m {\bf 9. Theorem.} $ { R}_i({\cal O}_V(H))$ occurs as  a subsheaf
of ${\cal O}_V(irH+iH+iK_V)$ generated by global sections.
\m As an incidental remark, it can be proven from the definitions that $R({\cal O}_V(H))$
is just the same as $R$ up to twisting, such that
$R_i({\cal O}_V(H))$ can be identified with ${\cal R}_i(i(r+1)H).$ 
\m Proof of Theorem 9.
Let us write the determinant in question in the form which resembles
 the differential in the Beilinson resolution  [9].
 For any divisor $E$ and
global sections $f_0,...,f_r$ of ${\cal O}_V(E)$ we can
construct a global section $w(f_0,...,f_r)$ of
${\cal O}_V( (r+1)E+K)$ by the formula
$$w(f_0,...,f_r)=\sum_{i=0}^r (-1)^i f_i\ df_0\wedge...\wedge
\widehat{(df_i)}\wedge...\wedge df_r$$
where the hat over the factor $df_i$ indicates that this
factor is to be deleted. We will now state and prove a lemma 
and then afterwards give the proof of Theorem 9.
\m {\bf 10. Lemma } The expression $w(f_0,...,f_r)$ is
 `homogeneous' in the sense that if $u$ is any
rational function
$$w(u f_0,...,u f_r)=u^{r+1}w(f_0,...,f_r).$$
\m Proof of Lemma 10.
The left side is the
right side plus
$$u^r\cdot \sum_{i=0}^r (-1)^if_i(\sum_{j=0}^{i-1} f_j df_1\wedge df_2 \wedge ..
.
\wedge du \wedge ... \wedge df_{i-1}\wedge df_{i+1}\wedge ... \wedge df_r
$$ $$+
\sum_{j=i+1}^rf_j df_1\wedge ...\wedge df_{i-1}
\wedge df_{i+1}\wedge ... \wedge du \wedge ... \wedge df_r).$$
Here the term $du$ replaces $df_j.$
\m There are two cancelling terms in which each pair of symbols
$df_i$ and $df_j$ is removed.
\m Proof of Theorem 9.
In  Lemma  10, observe
that  if one of the $f_i$ is taken to be $1$ then the expression
gives the wedge product of the remaining $f_i.$  
If the $f_i$ are allowed to range over
local sections of ${\cal O}$ the $w(f_0,...,f_r)$ span
a copy of the highest exterior power of the differentials of
$V$ modulo torsion. Whereas, if the $f_i$ range over
local sections of ${\cal O}_V(E)$ then the $w(f_0,...,f_r)$
generate a copy of the highest exterior power of the
differentials, twisted by $(r+1)E.$
\m If ${\cal F}$ is any subsheaf of
${\cal O}_V(E)$ and if the $f_i$ are allowed to
range over local sections but required to belong to the subsheaf ${\cal F},$
then the $w(f_0,...,f_r)$ generate a copy
of $F({\cal F})$ where $F$ is the functor we used above,
in the inductive definition of the $L_i.$
\m Thus we see
$${\cal F}\subset {\cal O}_V(E)\Rightarrow F({\cal F})\subset
{\cal O}_V((r+1)E+K)$$
Assuming inductively we have represented $$L_j({\cal O}_V(H))
\subset {\cal O}_V((r+2)^jrH+(r+2)^jH+(r+2)^jK_V)$$ for $j\le i.$
Consider the product sheaf
$$O_V(H)L_0({\cal O}_V(H))...L_i({\cal O}_V(H)).$$
This -- the  argument to which $F$ is to be applied  --
is a subsheaf of
$${\cal O}_V(H+(r+1)(1+(r+2)+(r+2)^2+...+(r+2)^i)H+(1+(r+2)+...+(r+2)^iK)$$
$$={\cal O}_V(H +((r+2)^{i+1}-1)H + {({(r+2)^{i+1}-1)}\over{r+1}}K).$$
We have seen  that $F$ multiplies both degrees by $(r+1)$ and adds $K,$
so
we can  represent $L_{i+1}$ in
$${\cal O}_V( (r+1)(r+2)^{i+1}H+(r+2)^{i+1}K)$$
as needed.
\m Now since ${ R}_i({\cal O}_V(H))$ is a product of the
$L_i$ for various values of $i$ we obtain the result desired
embedding.
\m Finally, now, we need to verify that the sections we
have used in the proof can be taken to be global.
Or, rather, that the global sections which we have
described generate the correct sheaves.
\m The issue is this. Choosing an affine open set $U\subset V$
we know the
$w(f_0,...,f_r)$ generate $F({\cal F})(U)$ if the $f_i$
range over the full module of sections ${\cal F}(U).$
But it is not necessarily true -- and easy to
find cases when it is false --   that the $w(f_0,...,f_r)$ generate
$F({\cal F})(U)$ if all we know is that the
$f_i$ range over a sequence of ${\cal O}_V(U)$-module
generators of ${\cal F}(U).$ In other words, we may
not merely state that the argument we have given works
when restricted to generating subsets.
\m However, we do know from [1] Proposition 11 that if we choose
a sequence of $k$ algebra generators of ${\cal O}_V(U),$
$x_0,x_1,...,x_n$ say,  such that $x_0=1,$ then $F({\cal F})(U)$ will be
generated by the $w(g_0,...,g_r)$ where the $g_i$
range over the pairwise products $x_if_j.$
\m Now, at any time when we are applying $F$
we apply it to a product ${\cal O}(H)L_0({\cal O}(H))...L_i({\cal O}(H))$
in which one of the factors is ${\cal O}(H).$
If we inductively assume all $L_i({\cal O}(H))$ are generated
by global sections we have specified already, and if we specify
$x_0,...,x_n$ to be a spanning sequence of global sections of ${\cal O}_V(H),$
then what we shall do is identify $V$ with the quasi-projective variety
one obtains using the embedding via the linear system coming
from $H.$ Thus  $V$ is covered by standard open subsets each of which
is isomorphic to a subset of affine space with the coordinates $x_0,...,x_n$
in which some $x_i$ is equal to one.
We shall choose such such a standard open set, and choose the numbering
so that $i=0.$ Now, the global sections of the product sheaf
which arise as $i+1$ fold products of the chosen sections on $U$ of the
separate sheaves
include all the products of the global sections of the $
L_j({\cal O}_V(H))$ which we have specified, ie times the section
$x_0=1$ of the leftmost sheaf, and also the same sections
multiplied by each of $x_1,...,x_n.$
\m Therefore upon applying the action of $F$ as a functor operating
on sheaves over this affine variety, the resulting sheaf $L_{i+1}(
{\cal O}_V(H))$
is indeed generated by the restrictions of the sections
$w(f_0,...,f_r)$ to our affine open subset,  where $f_0,...,f_r$ range
over the specified products of global sections.
\m This proves that the global sections produced by this process
restrict to module generators on each open subset of the standard
open cover of $V$.
\m And therefore that they generate the sheaves ${ R}_i({\cal O}_V(H))$
\m We also see, taking $T$ to be a fixed section of ${\cal O}_V((r+1)H),$
that the direct sum of the subsheaves of the ${\cal O}_V(irH+iH+K_V)$
is closed under multiplication, and in fact it is the same
sheaf of rings we have constructed already in earlier sections
on the affine open set where $T$ is not zero.

\m \centerline {11. Algebraic geometry proof of the theorem.}
\m Now we  give a quite easy algebraic geometry proof of 
theorem 3. We shall give only the important part
of the proof, in the most important case. Thus we assume $V$ is normal. 
We assume we 
are dealing with the codimension zero foliation only. These
are merely simplifying restrictions. And we give only
the implication which is most difficult, so we assume the chain of
Nash blowups  
$...\to V_{i+1}\to V_i \to ... \to V_0\to V$ stabilizes where
$V_0\to V$ is blowing up a sheaf of ideals $I,$ and we shall
prove the graded sheaf ${ R}$ is finite type.
\m 
Assume then that for some $i$ the map $ V_{i}\to V_{i-1}$ is
an isomorphism. Define divisors  $E, K_0, K_1,...$ on $V_i$    such
that ${\cal O}_{V_{i}}(E)$ is the pullback to $V_{i}$ of
$I$ (so $E$ is the negative of an exceptional Cartier divisor) and ${\cal O}_{V_{i}}(K_j)$ is the pullback to $V_{i}$
of the highest exterior power of the differentials of $V_j$  
mod torsion.
Note $K_j$ has been defined even for values of $j$ larger than $i$
since $V_i=V_{i+1}=V_{i+2}...$

\m 
Let $\pi$ be the composite map $V_{i}\to V$
We shall prove the following lemma by induction.
\m {\bf 11. Lemma.}
For all $j \ge 0$
$$\pi^*L_j(I)= {\cal O}_{V_{i}}(X_j)$$
where $X_j$ is the divisor
$$X_j=K_j + (r+1)[(r+2)^jE+\sum_{t=0}^{j-1} (r+2)^{j-1-t}K_t].$$
\m  Proof.
Write
$$L_j(I)=F(IL_0(I)L_1(I)..L_{j-1}(I)).$$
By  [1] corollary 3  
$$\pi^*F(IL_0(I)...L_{i-1}(I))\cong \pi^*(IL_0(I)...L_{i-1}(I))^{r+1}
\Lambda^r\Omega \cong {\cal O}_{V_{i}}(
(r+1)(E+\sum_{s=0}^{i-1}X_s)+K_i).$$ 
\m By the inductive hypothesis (the statement of the lemma
applied to
each  $X_s$) this equals ${\cal O}_V$ twisted by
$$(r+1)(E+\sum_{s=0}^{i-1}(K_s+ (r+1)[(r+2)^sE+\sum_{t=0}^{s-1}(r+2)^{s-1-t}K_t])+K_i
$$
$$=(r+1)(1 + (r+1) \sum_{s=0}^{i-1}(r+2)^s)E +(r+1)\sum_{s=0}^{i-1}K_s
+\sum_{t=0}^{i-1}\sum_{s=t+1}^{i-1}(r+1)^2(r+2)^{s-1-t}K_t + K_i$$
$$= (r+1)(1+(r+2)^i-1)E + (r+1)\sum_{t=0}^{i-1}((r+2)^{i-1-t}-1)K_t
+(r+1)\sum_{t=0}^{i-1}K_t+K_i$$
$$=(r+1)(r+2)^iE+(r+1)\sum_{t=0}^{i-1}(r+2)^{i-1-t}K_t+K_i$$
as required. QED
\m Now we resume the proof that $R(I)$ is finite type. Since $V_i=V_{i-1}$ then  $K_i=K_{i-1}.$
The formula just proven  with $K_i$ replaced by $K_{i-1}$ shows
$$X_i=K_{i-1}+(r+1)(r+2)^iE+(r+1)(r+2)^{i-1}K_0 + ... + (r+1)K_{i-1}.$$
Combining the terms $K_{i-1}$ and $(r+1)K_{i-1}$ we obtain
$$=(r+2)(K_{i-1}+ (r+1)(r+2)^{i-1}E+(r+1)(r+2)^{i-2}K_0+...+(r+1)K_{i-2}$$
$$=(r+2)X_{i-1}.$$

\m The graded sheaf
$$\pi^*{ R}(I)$$
in dimension $j$ when the base $r+2$ expansion of $j$ is
$j=a_0+a_1(r+2)+...+a_s(r+2)^s,$ is just
$${\cal O}_V( a_0X_0 + ... + a_sX_s).$$
Since $X_{j+1}=(r+2)X_j$ for large $j,$ the positive linear combinations
of all the divisors $X_0,X_1, ...,$ are spanned by
finitely many of the $X_i.$ In other words the $X_j$ generate 
a finitely generated monoid in the Weil divisor group. 
It follows that $\pi^*{ R}(I)$ is finite type over $k.$ 
The pushforward of each term
$$\pi_*\pi^*{ R}_j(I)$$
is the integral closure of ${ R}_j(I)$ by 
[11] page 219. This is because we are considering
the codimension zero foliation so $V_i$ is actually a nonsingular
variety for large $i.$ Actually [L] refers to sheaves of ideals
but the theorem applies to torsion free rank one coherent
sheaves, and the notion of integral closure extends.
Moreover $\pi_*\pi^*{ R}(I)$ is contained in 
the integral closure of the sheaf of rings $R(I).$  Applying
the result which states that an algebra over $k$
whose integral closure is finite type over $k$ is
itself finite type over $k,$  this shows that
the part of ${ R}(I) $ lying over each affine open set
is a graded $k$ algebra of finite type. 
This finishes the alternative,  algebraic geometry, proof of 
the more difficult implication of the theorem, subject to  the
inessential simplifying assumption that  $V$ is normal  and the
foliation is codimension zero. QED
\vfill\eject\centerline{12. A final result.}
\m We finish with a final result. Let $V$ be a normal irreducible
complex projective
variety with a resolvable foliation, and  $H$ a very ample divisor. 
Our hypothesis implies 
that there must exist a vector subspace $X\subset H^0(V, {\cal O}_V(iH))$
for some value of $i$ so that blowing up the base locus of $X$ resolves
the singularities of $V.$ We may view the elements of $X$
as homogeneous polynomials of degree $i.$ 
\m {\bf 12. Theorem.} There are functorial locally
closed conditions on the vector subspaces 
$X\subset |iH|,$ not vacuous for all $i,$ 
which when true, ensure
that blowing up the base locus of $X$ and one further Gaussian
blowup resolves the foliation.
\m {\bf Remark.} In case of the codimension zero foliation,
which is always resolvable (since singularities can be resolved)
the conditions ensure that the blowup of the base locus of $X$
is an immersed image of a smooth manifold. This is because
the map  $L_0^{r+2}\to L_1$ pulled back first to $V_0,$
the blowup of $I$, and then pulled back further to $V_1$, the Nash blowup of that,
yields up to twisting (the same on both sides), the map
from the pullback of the $r$'th exterior power of the differentials
of $V_0$ modulo torsion,  to the same for $V_1.$ This being
surjective implies that $V_0$ itself is an immersed image of
a smooth complex manifold.
\m {\bf Remark.} The reason for the word `functorial,' which
we will not define in this context precisely, 
is that, without some modification the theorem
is vacuously true. One needs  to rule out
proofs such as, to just arbitrarily choose a singleton resolving
vector space  $\{X\}$ for each
variety $V,$ using the axiom of choice. 

\m For the proof that the conditions are not
vacuous, though, we are allowed to make just
such an arbitary choice. Thus first imagine that we do start
with some choice of $X$, obtained in just that way, 
and let $I\subset {\cal O}_V(iH)$ be
the sheaf generated by $X.$ 
\m 
Now we merely reverse the argument, which we'll do in detail.
The sheaf $R(I)$ is finite type,
so  the subalgebra of $\Gamma(R({\cal O}_V(H))$ generated by
the particular global section generators we've described,
 whose degrees are powers
of $r+2,$ must be finite type;  because the full algebra of
global sections is integral over the subalgebra. Our
generating spaces, let us call 
them $X_j\subset H^0(V, {\cal O}_V((r+2)^j((r+1)iH+K)),$ because
they only occur in  power of $(r+2)$ degree, 
must satisfy $$X_j^{r+2}= X_{j+1}$$
for suitably large $j.$
For such a $j$ take now  $X=X_j$  and replace $i$ by $i(r+2)^j.$ We arrive
then at a situation where our new choice of $X$
satisfies  $X_0^{r+2}=X_1.$
\m
We now explicitly describe $X_0^{r+2}$ and $X_1$ in terms of the 
vector spaces  $X$ and $T=\Gamma(V, {\cal O}_V(H)).$
The former 
is the image of 
$$S^{r+2}\Lambda^{r+1}(X\otimes T )\to H^0(V, {\cal O}_V((r+2)(H+(r+1)iH+K_V))$$
sending a symmetric product of terms 
$(x_0\otimes t_0)\wedge ... \wedge (x_r\otimes t_{r})$
to the corresponding product of sections of the $w(x_0t_0,...,x_rt_r).$
The latter is 
the image
of the  map
$$\Lambda^{r+1}(X\otimes \Lambda^{r+1}(X\otimes T))\to H^0(V, {\cal O}_V(
(r+1)( H + (r+1)[H+(r+1)iH+K]+K)$$
with the same target,  sending 
  an exterior  product of monomials
$z_\alpha (x_0t_0\wedge ... \wedge x_rt_r)$ 
to the value of $w$ at the $z_\alpha w(x_0t_0,...,x_rt_r).$
\m The image of the first map is
contained in the image of the  second, even if we do not assume $X_0^{r+2}=X_1,$  and equality of the targets  implies
 by a determinantal rank formula in
terms of the coefficients of the homogeneous  polynomials. 
When it holds,  the map of sheaves generated by these global sections
 $L_0^{r+2}\to L_1$ is an isomorphism.
This map pulls back on $V_1$ (the Nash blowup of $V_0$ which is
the blowup of $I$), to the map of  the $r$'th exterior powers
of the differentials of $V_0$ pulled back to $V_1$ and reduced
modulo torsion to that of $V_1.$
This is surjective, and it follows that $V_0,$ the blowup 
of the base locus of $X,$ is 
an immersed image of a nonsingular variety.

\eject
\noindent [1] J. Moody, {\it On Resolving Singularities}, J. London Math Soc
(3) 64 (2001) 548-564 

\noindent
[2] Wahl, Jonathan M.(1-NC)
{\it The Jacobian algebra of a graded Gorenstein singularity}
Duke Math. J. 55 (1987), no. 4, 843--871 

\noindent 
[3] Cutkosky and Srinivas, {\it On a problem of Zariski
on dimensions of linerar systems.} Ann. of Math., 137 (1993) 531-559

\noindent
[4] O. Zariski, {\it The reduction of the singularities of an algebraic
surface}, Annals of Math 40 (1939) 699-689

\noindent
[5] L. Ein,{\it  The ramification divisors for branched coverings of} $  P\sp{n}\sb{k}$.
Math. Ann. 261 (1982), no. 4, 483--485.

\noindent
[6] S.  Kobayashi , T. Ochiai,  {\it    Characterizations of complex projective spaces and hyperquadrics}, J. Math. Kyoto Univ., Vol. 13-1 (1973), pp. 31-47.

\noindent
[7] K. Cho, Y.  Miyaoka, N. I.  Shepherd-Barron, {\it 
Characterizations of projective space and applications to complex symplectic manifolds. } 
Adv. Stud. Pure Math., 35 (2002) Math. Soc. Japan, Tokyo

\noindent
[8] Kebekus, Stefan {\it 
Characterizing the projective space after Cho, Miyaoka and Shepherd-Barron.} 
Complex geometry (Gottingen, 2000), 147-155, Springer, Berlin, 2002

\noindent
[7] Y. T. Siu, S-T Yau, {\it 
Compact Kahler manifolds of positive bisectional curvature.
} Invent. Math. 59 (1980), no. 2, 189--204

\noindent
[8] Baum, Fulton, Macphereson, {Riemann-Roch for singular
varieties} Publ. Math. IHES 45: 101-167, 1975

\noindent 
[9] A.A. Beilinson, {\it Coherent sheaves on Pn and problems in
linear algebra}, Funcsional Anal. i. prilozhen 12 (3) 1979 68-69

\noindent
[10] A. Beauville, J.-Y. Merindol, J.-Y.
{\it Sections hyperplanes des surfaces $K3$}
Duke Math. J. 55 (1987), no. 4, 873--878.

\noindent
[11] R. Lazarsfeld, positivity in algebraic geometry  II, Ergebnisse de Mathematik vol 48 (2004)

\noindent
[12]  H. Grauert O. Riemenschneider, {\it Verschwindungssatze fur analytische
Kohomologiegruppen  auf komplexen Raumen}, Invent Math 11 (1970) 263-292 

\noindent
[13] G. A. Gonzales-Springberg, {\it Resolution de Nash des points doubles rationnels} Annales de l'institut Fourier, 32 no. 2 (1982), p. 111-178

\noindent
[14] H. Hironaka, {\it On Nash blowing up},  Arithmetic and Geometry II,
Birkhauser 1983, 103-111

\noindent
[15]  M. Spivakovsky, {\it Sandwiched singularities and desingularization of surfaces by normalized Nash transformations.} Ann. of Math. (2) 131 (1990), no. 3, 411-491.

\noindent
[16 ]  F. A. Bogomolov, M. McQuillan,{\it Rational curves
on foliated varieties}  IHES preprint, February 2001

\noindent [17] Yu. Ilyashenko, {\it Centennial history of Hilbert's 16'th
problem} AMS Bulletin 39 (3) 301-354 (2002)

\noindent [18] M. Beltrametti, A. Sommese,  The adjunction theory of
complex projective varieties, DeGruyter, 1995
\end{document}